\documentclass[11pt]{amsart}
\usepackage{palatino, mathpazo}
\usepackage[mathscr]{eucal}
\usepackage{amsmath,amsthm,amsfonts,amssymb, mathrsfs, mathtools, bm}
\usepackage[all]{xy}
\usepackage{hyperref}
\usepackage{paralist}
\usepackage{microtype}

\newtheorem{theorem}{Theorem}[section]
\newtheorem{lemma}[theorem]{Lemma}
\newtheorem{proposition}[theorem]{Proposition}
\newtheorem{notation}[theorem]{Notation}
\newtheorem{corollary}[theorem]{Corollary}
\theoremstyle{remark}
\newtheorem{remark}[theorem]{Remark}

\numberwithin{equation}{section}

\newcommand{\depth}{\mathrm{depth}\,}
\newcommand{\rtop}{\mathrm{top}}
\newcommand{\Coker}{\mathrm{Coker}}
\newcommand{\Cl}{\mathrm{Cl}}
\newcommand{\Weil}{\mathrm{Weil}}
\newcommand{\Cart}{\mathrm{Cart}}
\newcommand{\rank}{\mathrm{rank}}
\newcommand{\Gr}{\mathrm{Gr}}

\newcommand{\Def}{\mathrm{Def}}

\DeclareMathOperator{\td}{td}

\newcommand{\hj}{\hat{j}}

\newcommand{\hZ}{\widehat{Z}}

\newcommand{\bC}{\mathbb{C}}

\newcommand{\bP}{\mathbb{P}}

\newcommand{\cO}{\mathcal{O}}

\newcommand{\sF}{\mathscr{F}}
\newcommand{\sV}{\mathscr{V}}
\newcommand{\sH}{\mathscr{H}}

\newcommand{\sZ}{\mathscr{Z}}

\title[On deformations of $1$-convex manifolds]{A remark on deformations of $1$-convex manifolds with exceptional curves}

\author[S.-S.\ Wang]{Sz-Sheng Wang}
\address[]{Institute of Mathematics, Academia Sinica, Taipei, Taiwan}
\email{sswangtw@gate.sinica.edu.tw}

\keywords{Deformation of singularities, $1$-convex manifolds, Small resolutions, Cohen--Macaulay singularities}

\subjclass[2020]{32S30, 14B07, 13H10, 32F10, 32S45, 32S35}

\begin{document}

\maketitle

\begin{abstract}
A formula for the dimension of the smoothing component of a $3$-dimensional isolated Cohen--Macaulay singularity is shown. We apply this formula for a $1$-convex threefold with a connected exceptional curve which is blown down to a terminal Gorenstein singularity.
\end{abstract}

\section{Introduction} \label{intro}

An $n$-dimensional complex manifold $X$ is called a \emph{$1$-convex} (or strongly pseudoconvex) $n$-fold if there is a proper surjective morphism $\pi$ from $X$ onto a Stein space $V$ with $\pi_{\ast} \cO_X \cong \cO_V$ and a finite subset $\Sigma \subseteq V$ such that $X \setminus E \to V \setminus \Sigma$ is biholomorphic where $E = \pi^{- 1} (\Sigma)$. We call $E$ the \emph{exceptional set} and $\pi \colon (X, E) \to (V, \Sigma)$ the \emph{Remmert reduction}.

In \cite[Theorem 2]{Laufer80}, Laufer proved that if $X$ is a $1$-convex manifold with a $1$-dimensional exceptional set $E$ then it has (the germ of) the miniversal deformation space $\Def (X)$. It is a natural question to ask for a formula for the dimension of $\Def (X)$.

We may assume without loss of generality that the curve $E$ is connected, i.e., $\Sigma = \{p\}$ and thus the germ $(V, p)$ is an isolated singularity. In the case of surfaces, Stephen Yau gave a formula of $\dim \Def (X)$ in \cite[Theorem 3.1]{YauS81} for an isolated hypersurface singularity $(V, p)$ and Wahl in \cite[Theorem 3.13 (c)]{Wahl81} for a smoothable, normal surface singularity. 

In the present note, we will generalize the results of \cite{Wahl81,YauS81}. The main result of this note is Theorem \ref{mainCM3} in Section \ref{pf_sec} for a $3$-dimensional isolated Cohen--Macaulay singularity which is smoothable. As an application, we derive a formula of $\dim \Def (X)$ for a certain $1$-convex threefold $X$:

\begin{theorem}\label{main1conv}
Let $X$ be a $1$-convex threefold with a connected exceptional curve $E$, and $\pi \colon (X, E) \to (V, p)$ the Remmert reduction as above. If $K_X$ is $\pi$-trivial, then
\begin{equation*}
    \dim \Def (X) = \tau - \frac{\mu + \sigma}{2}
\end{equation*}
where $\mu$ (resp.~$\tau$) is the Milnor (resp.~Tjurina) number of the singularity $(V, p)$ and $\sigma$ is the rank of the local divisor class group $\Cl (\cO_{V, p})$.
\end{theorem}

Here the canonical divisor $K_X$ is said to be  $\pi$-trivial if the intersection number of $K_X$ with every irreducible component of $E$ is zero. In such a case, the exceptional curve $E$ blows down to a terminal Gorenstein singularity, so to an isolated cDV (compound Du Val) singularity \cite{Reid83}. We remark that in general $(V, p)$ need not even be Cohen--Macaulay (see for example \cite[p.626 (1)-(4)]{Ando96} and \cite[Example 3.2]{Ando13}).

We will relate $\dim \Def (X)$ to Du Bois invariants $b^{r, s} (V, p)$ of the rational isolated hypersurface singularity $(V, p)$, which were introduced in \cite{St97} (see \S \ref{DB_subsec}).

\begin{corollary}\label{corodim}
Under the hypotheses of Theorem \ref{main1conv}, the dimension of $\Def (X)$ is equal to $b^{2, 1} (V, p)$. Moreover, the following are equivalent:
\begin{enumerate}
    \item\label{corodim1} The $(V, p)$ is an ordinary double point.
    \item\label{corodim2} $\dim  \Def(V) = \dim \Def(X) + 1$.
    \item\label{corodim3} $b^{1, 1} (V, p) = 0$.
    \item\label{corodim4} $b^{2, 1} (V, p) = 0$.
\end{enumerate}
\end{corollary}

The conditions \eqref{corodim2} and \eqref{corodim3} occur in the deformation theory of singular Calabi--Yau threefolds (see \cite[Lemma 3.6]{Gross97} and \cite[Theorem 2.2]{NS95}), and the condition \eqref{corodim4} was studied in \cite[Theorem 5.4, Remark 5.12]{St06}. There are similar results in \cite[\S 3]{FL22} (see Remark \ref{rmk_FL}).

We close this introduction with a few remarks on the techniques used in this paper. We will compare the Euler characteristics on the smoothing and the resolution (see Proposition \ref{lem-RRdef} and \ref{defect_dwR}) by using a globalization property of smoothings, proved by Looijenga \cite{Looij86}. Then we derive a formula of $\dim \Def (X)$ for a $1$-convex $n$-fold with a connected exceptional curve (cf.~Remark \ref{rmk_1conv}). It can be computed in principle by the Riemann--Roch theorem.

The Riemann--Roch defect for isolated singularities has been studied in \cite[\S 3]{Looij86}. However, the scissor relation (c) in the proof of \cite[Theorem 3.3]{Looij86} only holds for certain singularities (see Remark \ref{rmk_Loo}). Fortunately, the theorem of Riemann--Roch on threefolds is easy to compute (see Proposition \ref{RR}). Based on the method of Laufer, Wahl and Looijenga \cite{Looij86,Wahl81}, we will apply the Riemann--Roch formula to prove our results (see \S \ref{pf_sec}).

\medskip

\emph{Acknowledgements.}
The author is greatly indebted to Prof.~E.~Looijenga for answering the author's questions and Prof.~Chen-Yu Chi for many useful discussions. The author thanks Shing-Tung Yau Center of Southeast University and Institute of Mathematics at Academia Sinica for providing support and a stimulating environment, and was supported by the Fundamental Research Funds for the Central Universities 2242020R10048 and partially supported by the Ministry of Science and Technology (MOST) grant 110-2811-M-001-530.


\section{Preliminaries}

\subsection{Du Bois invariants}\label{DB_subsec}
We start by recalling the setup from \cite{St97}. Let $(Y, D) \to (V, p)$ be a good resolution of an isolated singularity $(V, p)$ of pure dimension $n$, i.e., $D$ is a divisor with simple normal crossings on $Y$. We fix a representative $Y \to V$ with $V$ a contractible Stein space.

Steenbrink defined invariants $b^{r, s}$ of $(V, p)$, called the \emph{Du Bois invariants} (see \cite[\S 2]{St97}), by
\begin{equation*}
    b^{r, s} (V, p) \coloneqq \dim H^s (Y, \Omega_Y^r (\log D) (- D))
\end{equation*}
for $r \geqslant 0$ and $s \geqslant 1$. These invariants do not depend on the choice of the good resolution, as they can be defined in terms of the filtered de Rham complex (\cite{DB81}, \cite[(3.5)]{St83}).

There are another invariants of $(V, p)$, defined in \cite[\S 3]{St97}. By \cite[Theorem 1.9]{St83}, there is a mixed Hodge structure on $H^k (V \setminus \{p\}) \cong H^{k + 1}_{\{p\}} (V)$. The Hodge filtration $F^{\cdot}$ on the local cohomology $H^{k + 1}_{\{p\}} (V)$ arises from a spectral sequence 
\begin{equation} \label{spE1}
    E_1^{p q} = H^q (D, \Omega_Y^p (\log D) \otimes \cO_D) \Rightarrow H^{p + q + 1}_{\{p\}} (V, \bC)
\end{equation}
which degenerates at the $E_1$-term. Therefore its Hodge filtration defines invariants $l^{r, s}$ of the singularity $(V, p)$ by $l^{r,s} \coloneqq \dim E_1^{s, r}$.

The following proposition follows from Lemma $2$, Theorem $4$ and $6$ in \cite{St97}. Note that the isolated rational singularity $(V, p)$ is Du Bois, i.e., $b^{0, s} = 0$ for $0 < s < n$ (see \cite[(3.7)]{St83}).

\begin{proposition}\label{prop_tm-DB}
Let $(V, p)$ be a 3-dimensional, rational isolated hypersurface singularity. Then
\begin{center}
    $\tau = b^{2, 1} + b^{1, 1} + l^{1, 1}$ and $\mu = 2 b^{1, 1} + l^{1, 1}$.
\end{center}
\end{proposition}

We remark that the main results in \cite{St97} express the Tjurina and Milnor numbers of certain isolated singularities in terms of Du Bois invariants $b^{r, s}$ and $l^{r, s}$. We state the proposition above considering such hypersurfaces for simplicity, which is all we need in the proof of Corollary \ref{corodim}.

Let $\sigma (V, p)$ denote the rank of the local divisor class group $\Cl (\cO_{V, p})$, i.e.,
\[
    \sigma (V, p) = \rank (\Weil (V, p) / \Cart (V, p)),
\]
where $\Weil (V, p)$ is the Abelian group of Weil (resp.~Cartier) divisors of the singularity $(V, p)$. It is a finite number if $(V, p)$ is a rational singularity (see \cite[Lemma 1.12]{Kawamata88}).

\begin{proposition}\label{prop_l11}
If $(V, p)$ is a 3-dimensional, rational isolated hypersurface singularity, then $\sigma(V, p) = l^{1, 1} (V, p)$.
\end{proposition}

\begin{proof}
Recall that $l^{1, 1} (V, p) = \dim \Gr_F^1 H^{3}_{\{p\}} (V, \bC)$ is defined by the spectral sequence \eqref{spE1}. Since $(V, p)$ is rational, we have $l^{0,i} = l^{i,0} = 0$ for all $i$ by \cite[Lemma 2]{St97} and in particular $E_{\infty}^{02} = E_{\infty}^{20} = 0$. Hence we get
$H^{3}_{\{p\}} (V, \bC) = \Gr_F^1 H^{3}_{\{p\}} (V, \bC)$. Then the proposition follows from the fact that $\sigma (V, p) = \dim H^2 (V \setminus \{p\}, \bC)$ (see \cite[(6.1)]{Flenner81} or the proof of \cite[Proposition 3.10]{NS95}).
\end{proof}

\begin{remark}\label{rmk_irrdef}
Let $\pi \colon X \to V$ be a resolution of a $3$-dimensional rational isolated singularity $(V, p)$. If the exceptional set $E$ of $\pi$ has dimension $1$, then the number of irreducible components of $E$ equals $\sigma (V, p)$ (cf.~\cite[Lemma 3.4]{Kawamata88} and \cite[Remark 2.8, 2.10]{SSW20}).
\end{remark}


\subsection{The Riemann--Roch defect}\label{RRdef_subsec}

Let $(V, p)$ be an isolated normal singularity of pure dimension $n \geqslant 2$ which is smoothable. 

Let $f \colon \sV \to \Delta \subseteq \bC$ be a good representative of a smoothing of $V = V_0$. According to a globalization theorem of Looijenga \cite[Appendix]{Looij86}, it follows that $f$ is a restriction of a projective flat family $F \colon \sZ \to \Delta$ which is smooth outside the point $\{p\}$. We wirte $Z$ for the fiber $Z_0$ and $Z_t$, $t \neq 0$, for non-singular fiber. Let $S$ be the smoothing component on which the smoothing $f \colon \sV \to \Delta$ takes place and $\beta_f \coloneqq \dim S$ (cf.~\cite[(4.1)]{Wahl81}).

The following lemma is known in \cite[(3.8)]{Wahl81}, and we recall the argument for the convenience of the reader. 

\begin{lemma} \label{lem-def}
With the above hypothesis and notation, for $t \neq 0$ we have 
\[
    \chi (\Theta_{Z}) = \chi (\Theta_{Z_t}) + \beta_f.
\]
\end{lemma}

\begin{proof}
By a conjecture of Wahl proved by Greuel and Looijenga \cite[(2.6)]{GL85}, we have
\begin{equation}\label{b_inv}
    \beta_f = \dim_{\bC} \Coker (\Theta_{\sZ / \Delta, p} \otimes \cO_{Z, p} \to \Theta_{Z, p}).
\end{equation}
Note that the natural morphism is injective because the sheaf of relative derivations $\Theta_{\sZ / \Delta}$ has depth $\geqslant 2$. Since $\Theta_{\sZ / \Delta}$ is flat over $\Delta$ and induces $\Theta_{Z_t}$ for $t \neq 0$, the lemma follows from the additivity of the Euler characteristic and the fact that the function $\chi (\Theta_{\sZ / \Delta} \otimes \cO_{Z_t})$ is constant on $t \in \Delta$ \cite[III, Theorem 4.12]{BS76}.
\end{proof}

\begin{remark}\label{rmk_Tju}
If $(V, p)$ is an isolated complete intersection singularity, then it is unobstructed. Hence $\beta_f$ is the Tjurina number $\tau$ of $(V, p)$, which is independent of the smoothing $f$.
\end{remark}

\begin{notation}\label{resl_nota}
Let $\pi \colon X \to V$ be any resolution of the singularity $(V, p)$. By gluing the resolution $\pi$ and the identity map of $Z \setminus \{p\}$, we get a resolution $\Pi \colon \hZ \to Z$.
\end{notation}

The following lemma will be useful in the sequel.

\begin{lemma}\label{LerayRes}
Let $\hZ$ and $\Pi$ be as in Notation \ref{resl_nota}. For a coherent sheaf $\sF$ on $\hZ$, we have 
\begin{equation*} 
    \chi (\sF) = \chi (\Pi_{\ast} \sF) + \sum_{1 \leqslant i \leqslant n - 1} (- 1)^i \dim (R^{i} \Pi_{\ast} \sF)_{p}.
\end{equation*}
\end{lemma}

\begin{proof}
Notice that $R^{i} \Pi_{\ast} \sF$ is supported on the point $\{p\}$ for $i > 0$. The lemma follows from the Leray spectral sequence for $\Pi$ and $\sF$. 
\end{proof}

We recall the notion of equivariant resolutions (in the sense of Hironaka). A resolution $\pi \colon X \to V$ is called \emph{equivariant} if the natural injection $\pi_{\ast} \Theta_X \hookrightarrow \Theta_V$ is an isomorphism. Note that such resolutions always exists (see \cite{Hironaka77} or \cite[p.14]{Kawamata85}).

\begin{proposition} \label{lem-RRdef}
With the above hypothesis and Notation \ref{resl_nota}, we assume that $\pi$ is equivariant. Then, for $t \neq 0$ we have
\begin{equation*}
    \chi (\Theta_{Z_t}) - \chi (\Theta_{\hZ}) = \sum_{1 \leqslant i \leqslant n - 1} (- 1)^{i - 1} h^{i} (\Theta_X) - \beta_f,
\end{equation*}
and $h^i (\Theta_X) = 0$ for all $i > r$ if the resolution $\pi$ has relative dimension $\leqslant r$.
\end{proposition}

Here and subsequently, $h^i (\sF)$ denotes the dimension of $H^i (X, \sF)$ for a coherent sheaf $\sF$ on $X$.

\begin{proof}
We apply Lemma \ref{LerayRes} to $\sF = \Theta_{\hZ}$, and replace $(R^{i} \Pi_{\ast} \Theta_{\hZ})_{p}$ by $h^i (\Theta_X)$ because $V$ is Stein. The first assertion follows from the fact that $\Pi$ is equivariant and Lemma \ref{lem-def} and the second assertion from the theorem on formal functions \cite[III, Theorem 3.1]{BS76}. 
\end{proof}

\begin{remark}\label{rmk_1conv}
Let $X$ be a $1$-convex $n$-fold with $1$-dimensional exceptional set. By \cite[Theorem 2]{Laufer80}, the miniversal deformation space $\Def (X)$ is smooth of dimension $h^1 (\Theta_X)$. It is also known that the Remmert reduction $\pi \colon X \to V$ is equivariant (cf.~\cite[(3.1)]{Friedman86}). If $V$ has a smoothing $f \colon \sV \to \Delta$, then by Proposition \ref{lem-RRdef} we get 
\[
    \dim \Def(X) = \chi (\Theta_{Z_t}) -  \chi (\Theta_{\hZ}) + \beta_f.
\]
For $n = 2$, this was given in \cite[(3.10.3)]{Wahl81}.
\end{remark}

Let $F$ be the Milnor fiber of the smoothing $f \colon \sV \to \Delta$, i.e., it is the fiber of the topological fiber bundle $f \colon \sV \setminus V \to \Delta \setminus \{0\}$. The middle Betti number $b_n (F)$ of $F$ of the smoothing $f$ is called the \emph{Milnor number} $\mu_f$ of $f$. If we write $E$ for the exceptional fiber $\Pi^{- 1} (p)$, then the difference of topological Euler characteristics $\chi_{\mathrm{top}}$ of $E$ and $F$ equals the global topological defect (see \cite[(3.5.3)]{Wahl81}),  
\begin{equation}\label{defect_top}
    \chi_{\mathrm{top}}(F) - \chi_{\mathrm{top}}(E) = \chi_{\mathrm{top}}(Z_t) - \chi_{\mathrm{top}} (\hZ).
\end{equation}

We remark that if the isolated singularity $(V, p)$ is a complete intersection then $\mu_f$ depends only on the singularity. In this case, we denote it by $\mu (V, p)$. Moreover, the Milnor fiber $F$ is $(n - 1)$-connected (see, e.g., \cite[(5.8)]{Looij84}) and there are only two non-vanishing Betti numbers $b_0 (F) = 1$ and $b_n (F)$, and thus
\begin{equation}\label{milnorNum}
        \chi_{\rtop} (F) = 1 + (- 1)^n \mu (V, p).
\end{equation}

\subsection{Dual of dualizing sheaves}

Keep the notation as Section \ref{RRdef_subsec}. We further assume that $(V, p)$ is Cohen--Macaulay. Then there is a relative dualizing sheaf $\omega_{\sZ / \Delta}$ for the globalization $\sZ \to \Delta$ of the smoothing $f$. It is flat over the $1$-dimensional disk $\Delta$. The dual $\omega_{\sZ / \Delta}^{\vee}$ is still torsion-free and thus flat over $\Delta$, and induces $\omega_{Z_t}^{\vee}$ for $t \neq 0$ and an inclusion $\omega_{\sZ / \Delta}^{\vee} \otimes \cO_Z \hookrightarrow \omega_{Z}^{\vee}$.

Wahl introduced the notion of \emph{$\omega^{\vee}$-constant deformations} \cite[(1.4)]{Wahl80}, to which the dual of the dualizing differentials lifts, and an invariant of the smoothing $f$ \cite[(3.7)]{Wahl81},
\begin{equation}\label{a_inv}
    \alpha_f \coloneqq \dim_{\bC} \Coker (\omega_{\sZ / \Delta, p}^{\vee} \otimes \cO_{Z, p} \hookrightarrow \omega_{Z, p}^{\vee}).
\end{equation}
Note that the smoothing $f$ is $\omega^{\vee}$-constant if and only if $\alpha_f = 0$ (automatic if $(V, p)$ is Gorenstein). By semicontinuity \cite[III, Theorem 4.12]{BS76}, we get for $t \neq 0$,  
\begin{equation}\label{defect_dwD}
    \chi (\omega_{Z_t}^{\vee}) = \chi (\omega_{\sZ / \Delta}^{\vee} \otimes \cO_Z) = \chi (\omega_Z^{\vee}) - \alpha_f.
\end{equation}

Next we treat the case of resolutions.

\begin{lemma}\label{dwR}
Under the above hypotheses, for the resolution $\Pi \colon \hZ \to Z$ as in Notation \ref{resl_nota}, we have a natural inclusion $\Pi_{\ast} (\omega_{\hZ}^{\vee}) \hookrightarrow \omega_{Z}^{\vee}$. Furthermore, it is an isomorphism if the exception set of $\Pi$ has codimension $\geqslant 2$.

\end{lemma}

\begin{proof}
Let $E = \Pi^{- 1} (p)$, and let $\hj \colon \hZ \setminus E \hookrightarrow \hZ$ be the natural inclusion. Consider the exact sequence of local cohomology sheaves \cite[Corollary 1.9]{Hartshorne67}
\begin{equation}\label{wdulLoC}
    0 \to \sH_E^0 (\omega_{\hZ}^{\vee}) \to \omega_{\hZ}^{\vee} \to \hj_{\ast} (\hj^{\ast} \omega_{\hZ}^{\vee}) \to \sH_E^1 (\omega_{\hZ}^{\vee}) \to 0
\end{equation}
and note that $\sH_E^0 (\omega_{\hZ}^{\vee}) = 0$ by a depth argument (cf.~\cite[Theorem 3.8]{Hartshorne67}). Furthermore, $\sH_E^1 (\omega_{\hZ}^{\vee}) = 0$ if $E$ has codimension $\geqslant 2$.

By \eqref{wdulLoC}, we have  $\Pi_{\ast} (\omega_{\hZ}^{\vee}) \hookrightarrow \Pi_{\ast} \hj_{\ast} (\hj^{\ast} \omega_{\hZ}^{\vee})$. Let $U \coloneqq Z \setminus \{p\}$ be the smooth locus and $j \colon U \hookrightarrow Z$ the inclusion. Then, since $\depth \omega^{\vee}_{Z, p} \geqslant 2$, we get $\omega_{Z}^{\vee} \xrightarrow{\sim} j_{\ast} (\omega_U^{\vee})$ by $\sH_{\{p\}}^i (\omega_{Z}^{\vee}) = 0$ for $i =0, 1$. Thus the proposition follows from $\Pi_{\ast} \hj_{\ast} (\hj^{\ast} \omega_{\hZ}^{\vee}) = j_{\ast} \Pi_{\ast}  (\hj^{\ast} \omega_{\hZ}^{\vee}) = j_{\ast} (\omega_U^{\vee})$.
\end{proof}

\begin{remark}
If $(V, p)$ has dimension two and $\Pi \colon \hZ \to Z$ is the minimal resolution, then one always has $\Pi_{\ast} (\omega_{\hZ}^{\vee}) \xrightarrow{\sim} \omega_{Z}^{\vee}$ by \cite[(3.5)]{Wahl80}.
\end{remark}

\begin{proposition}\label{defect_dwR}
Let $X$, $\hZ$, $\pi$ and $\Pi$ be as in Notation \ref{resl_nota}. Then we have 
\[
    \chi (\omega_{Z_t}^{\vee}) - \chi (\omega^{\vee}_{\hZ}) = \sum_{1 \leqslant i \leqslant n - 1} (- 1)^{i - 1} h^i (\omega_X^{\vee}) - \alpha_f + \gamma_{\pi}
\]
with $\gamma_{\pi} = \dim_{\bC} \Coker (\Pi_{\ast} (\omega_{\hZ}^{\vee})_{p} \hookrightarrow \omega_{Z, p}^{\vee})$. 
\end{proposition}

\begin{proof}
By Lemma \ref{dwR} and the additivity of the Euler characteristic, we get $\chi (\omega_Z^{\vee}) = \chi (\Pi_{\ast} \omega_{\hZ}^{\vee}) + \gamma_{\pi}$. Then the corollary follows from Lemma \ref{LerayRes}, \eqref{defect_dwD} and that $V$ is Stein.
\end{proof}

\begin{remark}\label{rmk_dw0}
Suppose that $(V, p)$ is a rational Gorenstein singularity, i.e., it is a canonical singularity of index $1$ (cf.~\cite[(6.4), (6.8)]{CKM88} or \cite[(3.3)]{Gross97}). In the case $\omega_X \cong \pi^{\ast} \omega_V$, we have $h^i (\omega_X^{\vee}) = 0$ for $i > 0$. Such $\pi$ is called a \emph{crepant resolution}.

Indeed, we may assume that $\omega_V \cong \cO_V$ for the Gorenstein singularity $(V, p)$. Then $(R^i \pi_{\ast} \omega_X^{\vee})_p \cong (R^i \pi_{\ast} \cO_X)_p = 0$ for $i > 0$ since $(V, p)$ is rational.
\end{remark}


\section{Main result} \label{pf_sec}
Suppose from now on that $(V, p)$ is a $3$-dimensional isolated normal Cohen--Macaulay singularity with a smoothing $f \colon \sV \to \Delta$. Given a resolution $\pi \colon X \to V$, one can define the \emph{geometric genus} $p_g (V, p)$ of $(V, p)$ by $h^2 (\cO_X)$. It is well known that the number is independent of the resolution. 

In the previous section, we have seen the invariants $\alpha_f$ and $\beta_f$ associated the smoothing $f$, defined as in \eqref{a_inv} and \eqref{b_inv} respectively. If the resolution $\pi$ is equivariant, i.e., $\pi_{\ast} \Theta_X \cong \Theta_V$, then the following main result are going to relate these invariants $\alpha_f$ and $\beta_f$ with certain numbers induced by $\pi$. This is a generalization of \cite[Theorem 3.13 (c)]{Wahl81}.

\begin{theorem} \label{mainCM3}
With the above hypothesis and notation, we assume that $\pi$ is equivariant. Then we have
\begin{align*}
    \beta_f - \alpha_f = &h^1 (\Theta_{X}) - h^2 (\Theta_{X})  - 22 p_g (V, p) - \frac{1}{2} (\chi_{\mathrm{top}}(F) - \chi_{\mathrm{top}}(E)) \\
    &- (h^1 (\omega_X^{\vee}) - h^2 (\omega_X^{\vee})) - \gamma_{\pi}
\end{align*}
where $E$ is the exceptional set of the resolution $\pi$, $F$ the Milnor fiber of the smoothing $f$ and $\gamma_{\pi} = h^0 (\omega_{V}^{\vee} / \pi_{\ast} (\omega_{X}^{\vee}))$. 
\end{theorem}

Before stating the result to be proved, we need the following proposition, which follows from Hirzebruch--Riemann--Roch theorem and the observation that $\omega_{M}^{\vee} = \det \Theta_M$ and $\td_3 (\Theta_M) = (1 / 24) c_1 (\Theta_M) c_2 (\Theta_M)$ (cf.~\cite[Example 15.2.5]{Fulton98}).

\begin{proposition} \label{RR}
Let $M$ be a $3$-dimensional compact manifold. Then we have
\[
    \chi (\Theta_{M}) = \chi (\omega^{\vee}_{M}) - 22 \chi (\cO_M) + (1 /2) \chi_{\mathrm{top}} (M).
\]
\end{proposition}

We now obtain our main result:

\begin{proof}[Proof of Theorem \ref{mainCM3}]
We shall use the same notation as in Section \ref{RRdef_subsec}. By Proposition \ref{lem-RRdef}, for $t \neq 0$, we have 
\[
    \chi (\Theta_{Z_t}) - \chi (\Theta_{\hZ}) = h^1 (\Theta_{X}) - h^2 (\Theta_{X}) - \beta_f.
\]
Applying Proposition \ref{RR} to the left hand side gives
\begin{align*}
    \chi (\Theta_{Z_t}) -  \chi (\Theta_{\hZ}) =& \chi (\omega^{\vee}_{Z_t}) - \chi (\omega^{\vee}_{\hZ}) - 22 (\chi (\cO_{Z_t}) - \chi (\cO_{\hZ}))  \\
    &  + (1 / 2) (\chi_{\mathrm{top}} (Z_t) - \chi_{\mathrm{top}} (\hZ)).
\end{align*}
By \eqref{defect_top} and Proposition \ref{defect_dwR}, it suffices to show that
\begin{equation}\label{defect_O}
    \chi (\cO_{Z_t}) - \chi (\cO_{\hZ}) = - p_g (V, p).
\end{equation}
To do this, we notice that only $h^2 (\cO_X)$ can be non-zero for the normal Cohen--Macaulay singularity $(V, p)$, and $\chi (\cO_{Z_t}) = \chi (\cO_{Z})$ by semicontinuity \cite[III, Theorem 4.12]{BS76}. Then the equality \eqref{defect_O} follows from Lemma \ref{LerayRes}, and the proof is complete. 
\end{proof}

\begin{remark}\label{rmk_Loo}
Let the hypotheses be as in Theorem \ref{mainCM3} and $S$ the smoothing component on which the smoothing
$f$ takes place. Our result shows that
\[
    \dim S + (1 / 2) \chi_{top} (F) - \alpha_f
\]
is independent of $S$ (and depend only on $\pi \colon X \to V$). In particular, if one smoothing $f$ is $\omega^{\vee}$-constant (i.e., $\alpha_f = 0$), so are all smoothings on the same irreducible component $S$ (cf.~\cite[(3.14.2)]{Wahl81}).

As a consequence, the formula in \cite[4.4]{Looij86} is not correct in general.
\end{remark}

Theorem \ref{mainCM3} allows us to treat the case of certain $1$-convex threefolds.

\begin{proof}[Proof of Theorem \ref{main1conv}]
Recall that the exceptional set $E$ of the Remmert reduction $\pi \colon X \to V$ has dimension $1$ and $K_X$ is $\pi$-trivial. Then $(V, p)$ is a $3$-dimensional Gorenstein terminal singularity (see \cite{Reid83} or \cite[(16.2)]{CKM88}). In particular, it is a rational hypersurface singularity, and $p_g (V, p) = 0$ (cf.~\cite[Corollary 4.2]{Karras86}). According to \eqref{a_inv}, Lemma \ref{dwR} and Remark \ref{rmk_Tju}, it follows that $\alpha_f = \gamma_{\pi} = 0$ and $\beta_f$ equals the Tjurina number $\tau$ of $(V, p)$ for any fixed smoothing $f$.

We note that $\pi$ is crepant since the exceptional set $E$ contains no divisors. By Remark \ref{rmk_1conv}, \ref{rmk_dw0} and Theorem \ref{mainCM3}, we find that
\begin{align*}
    \dim \Def(X) = \tau + (1 / 2) (\chi_{\mathrm{top}}(F) - \chi_{\mathrm{top}}(E)). 
\end{align*}

It remains to prove that $\chi_{\mathrm{top}}(F) - \chi_{\mathrm{top}}(E)$ is equal to $- (\mu + \sigma)$ where $\mu$ is the Milnor number of $(V, p)$ and $\sigma$ is the rank of $\Cl (\cO_{V, p})$. Indeed, the irreducible components of $E$ are smooth rational curves, meeting transversally with no cycles \cite[Proposition 1]{Pinkham83}. By using the Mayer--Vietoris sequence, we get the number of the irreducible components of $E$ equals $\chi_{\mathrm{top}}(E) - 1$. Hence the desired formula follows from \eqref{milnorNum} and Remark \ref{rmk_irrdef}.
\end{proof}

\begin{proof}[Proof of Corollary \ref{corodim}]
First, under the hypotheses of Theorem \ref{main1conv}, we have seen that $(V, p)$ is an isolated cDV hypersurface singularity, and also have $\sigma > 0$ by Remark \ref{rmk_irrdef}. Observe that it is an ordinary double point if and only if $\mu = 1$; in particular, $\tau = \sigma = 1$.

By Theorem \ref{main1conv} and Proposition \ref{prop_tm-DB} and \ref{prop_l11}, we get that the dimension of $\Def (X)$ equals the Du bois invariant $b^{2, 1}$, and the condition \eqref{corodim1} clearly implies \eqref{corodim2}, \eqref{corodim3} and \eqref{corodim4}. 

Conversely, we can rewrite \eqref{corodim2} as $\mu + \sigma = 2$ and thus this implies \eqref{corodim1}. The condition \eqref{corodim3} implies \eqref{corodim1} by \cite[Theorem 2.2]{NS95} (see also \cite[p.1374]{St97}). 

Now suppose that $(V, p)$ is not an ordinary double point. Notice that $V$ can be considered as the total space of a deformation of a Du Val surface singularity. One can find a nontrivial small deformation of $X$ under which $E$ splits up into a finite disjoint union of smooth copies of $\bP^1$ with normal bundle $\cO_{\bP^1} (- 1)^{ \oplus 2}$ (see \cite[Proposition 1.1]{Wilson97} or \cite[p.679]{Friedman86}). Hence we get $\dim \Def (X) > 0$ and complete the proof of the equivalence of \eqref{corodim1} and \eqref{corodim4}. 
\end{proof}

\begin{remark}\label{rmk_FL}
In \cite[\S 3]{FL22}, Friedman and Laza proved similar results as Theorem \ref{main1conv} and Corollary \ref{corodim} with a different approach. 
\end{remark}

\bibliographystyle{plain}

\end{document}